\documentclass[11pt]{article}

\usepackage[margin=2.7cm]{geometry}
\usepackage[utf8]{inputenc}
\usepackage[T1]{fontenc}
\usepackage{lmodern}
\usepackage{amsmath,amssymb,amsthm,amsfonts,mathtools}
\usepackage{enumitem}
\usepackage{needspace}
\usepackage{aliascnt}
\usepackage{hyperref}
\usepackage[nameinlink,capitalize]{cleveref}
\usepackage{microtype}

\hypersetup{
    colorlinks=true,
    linkcolor=blue,
    citecolor=blue,
    urlcolor=blue,
    pdftitle={Generic local identifiability for ODE inverse problems from discrete observations},
    pdfauthor={Marian Petrica and Ionel Popescu}
}

\newtheorem{theorem}{Theorem}[section]
\newaliascnt{proposition}{theorem}
\newtheorem{proposition}[proposition]{Proposition}
\aliascntresetthe{proposition}
\newaliascnt{corollary}{theorem}
\newtheorem{corollary}[corollary]{Corollary}
\aliascntresetthe{corollary}
\theoremstyle{definition}
\newaliascnt{example}{theorem}
\newtheorem{example}[example]{Example}
\aliascntresetthe{example}

\newcommand{\R}{\mathbb{R}}

\newcommand{\Z}{\mathbb{Z}}
\newcommand{\Q}{\mathbb{Q}}
\newcommand{\rank}{\operatorname{rank}}

\title{Generic local identifiability for ODE inverse problems from discrete observations}

\author{
  Marian Petrica$^{1,2}$\\ 
  \small $^1$Gheorghe Mihoc--Caius Iacob Institute of Mathematical Statistics\\[-1mm]
  \small and Applied Mathematics of the Romanian Academy, Bucharest, Romania\\
  \small $^2$SEARCH-Vasc, Center of Excellence, Bucharest, Romania
  \and
  Ionel Popescu$^{3,4}$\\
  \small $^3$Faculty of Mathematics and Computer Science, University of Bucharest, Romania\\[-1mm]
  \small $^4$Simion Stoilow Institute of Mathematics of the Romanian Academy, Bucharest, Romania\\
}

\date{}

\begin{document}
\maketitle

\begin{abstract}
We study local parameter identification from finitely many measurements of an analytic output along a single ODE trajectory. At a fixed initial state and parameter, the time derivatives of the output sensitivity define a decreasing sequence of parameter kernels. If the sequence reaches zero, any prescribed pattern of distinct positive sampling times gives full parameter rank after a sufficiently small common rescaling, provided the number of times is at least the number of strict kernel drops. The orders of these drops determine the asymptotic orders of all singular values of the sampled Jacobian. We also give a sharp worst-case bound on the number of freely chosen
observation times needed for the sampled Jacobian to have full column
rank. At a generic point, the maximum parameter rank determines how
many local parameter coordinates can vary without changing the
observations. A quantitative inverse estimate bounds the error of every estimate in the certified neighborhood that fits the data at least as well as the true parameter. For a \(d\)-dimensional homogeneous linear system with known initial state and \(d\) equally spaced positive-time full-state observations, cyclicity for the sampled propagator is equivalent to local identifiability of the generator. In this case, one additional exact full-state sample at a time incommensurate with the grid spacing gives global recovery.
\end{abstract}



\section{Introduction}

A continuous output trajectory may determine the parameters of an ordinary differential equation locally. In practice, observations are taken at finitely many times, often on a uniform grid. We study which choices of observation times preserve local identifiability and how their spacing affects the conditioning of the inverse problem.

For an analytic model at a fixed initial state and parameter, we form a matrix by placing the time derivatives of the output sensitivity at \(t=0\) in successive row blocks. Suppose this matrix attains full column rank after finitely many derivatives and let \(s\) be the number of derivative orders at which its rank increases. For any \(0<c_1<\cdots<c_m\) with \(m\ge s\), observations at \(hc_1,\ldots,hc_m\) give a Jacobian of full column rank for all sufficiently small \(h>0\). In particular, the \(s\) equally spaced observations at \(h,2h,\ldots,sh\) suffice.

The increase in rank at derivative order \(k\) equals the number of
singular values of order \(h^k\) as \(h\downarrow0\). In particular,
if the last rank increase occurs at order \(k_*\), the smallest
singular value is of order \(h^{k_*}\), and the worst-case local
amplification of observation errors is of order \(h^{-k_*}\).

Under the same rank assumption, when the observation times can vary independently, almost every choice of \(p-r_0+1\) times yields a Jacobian of full column rank. Here \(p\) is the number of parameters and \(r_0\) is the maximum rank of the output sensitivity over the observation interval.

The use of time derivatives in identifiability builds on the power-series method of Pohjanpalo \cite{pohjanpalo1978system} and the observability rank condition of Hermann and Krener \cite{hermann1977nonlinear}, applied to parameters treated as constant states \cite{villaverde2016structural}. For freely chosen observation times, the closest result is Sontag's local theorem \cite[Sec.~2.2]{sontag2002differential}. It shows that at most \(p\) experiments suffice to preserve local distinguishability near a parameter where the maximal response rank is locally constant, even when that rank is below \(p\). His global \(2p+1\) theorem gives a single choice of experiments that separates all distinguishable parameter pairs. Van Willigenburg, Stigter and Molenaar \cite{vanwilligenburg2022sensitivity} compare rank conditions based on sampled sensitivities and Lie derivatives. The observer literature has also studied preservation of observability under fine sampling and the ill-conditioning that can accompany it \cite{ammar2004preservation,moraal1995observer}.

For homogeneous linear systems with full-state observations on a uniform grid, a cyclic orbit determines the sampled propagator and its real matrix logarithms describe the compatible generators \cite{stanhope2014identifiability,culver1966existence}. We show that cyclicity is equivalent both to invertibility of the sampled Jacobian and to local identifiability of the generator. Under cyclicity, one additional exact full-state observation at a time incommensurate with the grid spacing determines the generator globally.

\section{Finite observations and the sampling theorem}
\label{sec:sampling}

Let \(U\subset\R^d\), \(\Omega\subset U\times\R^p\) be open and let
\(f:\Omega\to\R^d\) and \(g:U\to\R^q\) be real analytic. We consider the
initial value problem and its observed output
\begin{equation}\label{eq:ode}
    \dot X(t)=f(X(t),\alpha),\qquad
    X(0)=x,\qquad
    y(t;\alpha,x)=g(X(t;\alpha,x)).
\end{equation}
The initial state \(x\) is known and is considered fixed when identifying the parameter vector \(\alpha\). Full-state observations correspond to \(q=d\) and \(g(z)=z\) for \(z\in U\).

Fix \(T>0\) and a reference pair \((\bar x,\bar\alpha)\in\Omega\) such that
\((X(t;\bar\alpha,\bar x),\bar\alpha)\) is defined and belongs to
\(\Omega\) for every \(t\in[0,T]\). By adjoining the equation
\(\dot\alpha=0\), we obtain the analytic vector field \((f,0)\) on
\(\Omega\). Its maximal flow is jointly analytic on its open flow domain
\cite[Chap.~I]{ilyashenko2008lectures}, which contains the compact set
\([0,T]\times\{(\bar x,\bar\alpha)\}\). Hence there are open balls
\(V_x\) and \(V_\alpha\) centered at \(\bar x\) and \(\bar\alpha\),
respectively, and an \(\epsilon_0>0\) such that
\(X(t;\alpha,x)\) exists and is jointly analytic on
\[
    (-\epsilon_0,T+\epsilon_0)\times V_x\times V_\alpha,
\]
with coordinates \((t,x,\alpha)\). In particular, \([0,T]\) is a common
interval of existence for all \((x,\alpha)\in V_x\times V_\alpha\).

For \(m\) observations, let
\[
    \mathcal T_m(T)
    =\{(t_1,\ldots,t_m)\in\R^m:0<t_1<\cdots<t_m<T\}.
\]
Given \(\tau\in\mathcal T_m(T)\) and \(x\in V_x\), the observation map
\(\Phi_{\tau,x}:V_\alpha\to\R^{mq}\) and its parameter Jacobian are
\begin{equation}\label{eq:observation-map}
    \Phi_{\tau,x}(\alpha)
    =\bigl(y(t_1;\alpha,x),\ldots,y(t_m;\alpha,x)\bigr),\qquad
    S(\tau,x,\alpha)=D_\alpha\Phi_{\tau,x}(\alpha).
\end{equation}
A parameter \(\alpha_0\) is locally identifiable for \((\tau,x)\) if
\(\Phi_{\tau,x}(\alpha)=\Phi_{\tau,x}(\alpha_0)\) implies
\(\alpha=\alpha_0\) for all \(\alpha\) sufficiently close to
\(\alpha_0\). If \(S(\tau,x,\alpha_0)\) has full column rank, we can
select \(p\) scalar components of \(\Phi_{\tau,x}\) whose Jacobian at
\(\alpha_0\) is invertible. The inverse function theorem applied to
these components gives local injectivity and, after restricting to a
smaller neighborhood, a Lipschitz inverse on the local image of
\(\Phi_{\tau,x}\). Local identifiability can also hold at a
rank-deficient point. For example, in the scalar model
\(\dot X=\alpha^3\), the observation \(X(t)=x+t\alpha^3\) at any \(t>0\)
determines \(\alpha\) uniquely, even at \(\alpha=0\), where its parameter
derivative vanishes. The sampling results below concern full column rank.

Fix \((x_0,\alpha_0)\in V_x\times V_\alpha\) and let
\(Z(t)=D_\alpha y(t;\alpha_0,x_0)\) be the output sensitivity at this
pair. The sampled Jacobian consists of the corresponding row blocks,
\[
    S(\tau,x_0,\alpha_0)
    =
    \begin{pmatrix}
        Z(t_1)\\
        \vdots\\
        Z(t_m)
    \end{pmatrix}.
\]
To express the time derivatives of \(Z\) at zero in terms of \(f\) and
\(g\), define
\[
    G_0(z,\alpha)=g(z),\qquad
    G_{k+1}(z,\alpha)=D_zG_k(z,\alpha)f(z,\alpha),\qquad k\ge0.
\]
Repeated differentiation along the solution gives
\[
    \frac{d^k}{dt^k}y(t;\alpha,x)
    =G_k(X(t;\alpha,x),\alpha).
\]
Joint analyticity allows us to interchange time and parameter
derivatives. Since \(X(0;\alpha,x_0)=x_0\) for every \(\alpha\), we obtain
\begin{equation}\label{eq:jet-kernels}
    \begin{gathered}
        Z(0)=0,\qquad
        Z_k:=Z^{(k)}(0)=D_\alpha G_k(x_0,\alpha_0),\quad k\ge1,\\
        K_0=\R^p,\qquad
        K_k=\bigcap_{j=1}^{k}\ker Z_j,\quad k\ge1.
    \end{gathered}
\end{equation}
Thus \(K_k\) consists of the parameter directions \(v\) for which the
first \(k\) time derivatives of \(Z(t)v\) vanish at zero. In particular,
\(Z_1=Dg(x_0)D_\alpha f(x_0,\alpha_0)\). Since
\(Z(t)/t\to Z_1\) as \(t\to0\), full column rank of \(Z_1\) guarantees
full column rank of \(Z(t)\) for every sufficiently small \(t>0\).
In this case, one observation at such a time suffices for local
identifiability.

\begin{proposition}[Time-derivative criterion and freely chosen times]
\label{prop:lie-witness}
There exists a finite set of observation times in \((0,T)\) for which
\(S(\tau,x_0,\alpha_0)\) has full column rank if and only if
\(\bigcap_{k\ge1}K_k=\{0\}\). Equivalently, \(K_k=\{0\}\) for some finite
\(k\). Under this condition, let
\[
    r_0=\max_{0<t<T}\rank Z(t).
\]
For every \(m\ge p-r_0+1\), the time tuples giving full column rank form
an open, dense subset of \(\mathcal T_m(T)\) of full Lebesgue measure.
\end{proposition}

\begin{proof}
The sequence \(K_0\supseteq K_1\supseteq\cdots\) can have only finitely
many strict inclusions, since each strict inclusion lowers the
dimension. It therefore eventually stabilizes at
\(\bigcap_{k\ge1}K_k\), which proves the equivalence of the two
conditions on the kernels.

If \(v\in\bigcap_{k\ge1}K_k\), all derivatives of \(Z(t)v\) at zero
vanish, including the value at zero. By analyticity, \(Z(t)v=0\) for
all \(t\in(0,T)\). A nonzero vector in the intersection therefore lies
in the kernel of every sampled Jacobian.

Now suppose that \(\bigcap_{k\ge1}K_k=\{0\}\). Choose
\(t_1\in(0,T)\) such that \(\rank Z(t_1)=r_0\). After choosing
\(t_1<\cdots<t_j\), write
\[
    L_j=\bigcap_{i=1}^{j}\ker Z(t_i)
\]
for the kernel of the Jacobian formed from these observations. If
\(L_j\) contains a nonzero vector \(v\), there is a time
\(t_{j+1}\in(t_j,T)\) with \(Z(t_{j+1})v\ne0\). Indeed, if \(Z(t)v\)
vanished throughout \((t_j,T)\), the analytic identity theorem would
force all its derivatives at zero to vanish, placing \(v\) in
\(\bigcap_{k\ge1}K_k\). The new time gives
\(L_{j+1}=L_j\cap\ker Z(t_{j+1})\subsetneq L_j\). Since
\(\dim L_1=p-r_0\), at most \(p-r_0\) additional times are needed to
reduce the common kernel to zero. Any further observations preserve
full column rank.

For fixed \(x_0,\alpha_0\), the Gram determinant
\(\det(S^{\mathsf T}S)\) is real analytic as a function of \(\tau\) on
the connected open simplex \(\mathcal T_m(T)\). It vanishes exactly
when \(S\) has deficient column rank. For each \(m\ge p-r_0+1\), the
construction gives a tuple at which this determinant is positive.
Its zero set is therefore relatively closed and has Lebesgue measure
zero \cite{mityagin2020zero}. The full-rank tuples form its open, dense
complement.
\end{proof}

\begin{theorem}[Prescribed sampling patterns and singular-value orders]
\label{thm:prescribed-grid}
Suppose that \(K_k=\{0\}\) for some \(k\), where \(K_k\) is defined in
\eqref{eq:jet-kernels}. Let \(k_1<\cdots<k_s\) be the derivative orders
at which \(K_{k_j}\subsetneq K_{k_j-1}\) and set
\[
    d_j=\dim K_{k_j-1}-\dim K_{k_j},
    \qquad j=1,\ldots,s.
\]
Thus \(d_j\) is the increase in rank when \(Z_{k_j}\) is added to the
preceding derivative blocks, and \(\sum_{j=1}^s d_j=p\).

Fix \(m\ge s\) and \(0<c_1<\cdots<c_m\). For \(0<h<T/c_m\), take
observations at the times \(\tau_h=(hc_1,\ldots,hc_m)\) and write
\(S_h=S(\tau_h,x_0,\alpha_0)\). Let \(\nu_1\le\cdots\le\nu_p\) be
the list obtained by repeating each \(k_j\) exactly \(d_j\) times.
With the singular values of \(S_h\) listed in decreasing order, there
exist \(a,b>0\) and \(0<h_0<\min\{1,T/c_m\}\) such that
\begin{equation}\label{eq:singular-orders}
    a h^{\nu_\ell}\le\sigma_\ell(S_h)\le b h^{\nu_\ell},
    \qquad \ell=1,\ldots,p,\qquad 0<h<h_0.
\end{equation}
In particular, \(S_h\) has full column rank for all sufficiently small
\(h>0\), with \(d_j\) singular values of order \(h^{k_j}\) for each
\(j\). Its smallest singular value is therefore of order \(h^{k_s}\).
The set of \(h\in(0,T/c_m)\) for which \(S_h\) has deficient column
rank is finite.
\end{theorem}

\begin{proof}
At each order \(k_j\), let
\[
    E_j=K_{k_j-1}\cap K_{k_j}^{\perp}.
\]
Then \(K_{k_j-1}=E_j\oplus K_{k_j}\) and \(\dim E_j=d_j\).
The kernels are constant between successive strict decreases, so
these splittings give the orthogonal decomposition
\[
    \R^p=E_1\oplus\cdots\oplus E_s.
\]
Since \(K_{k_j}=K_{k_j-1}\cap\ker Z_{k_j}\), the restriction of
\(Z_{k_j}\) to \(E_j\) is injective. In particular, \(d_j\le q\) and
hence \(p\le sq\le mq\).

For \(v_j\in E_j\), we have \(Z(0)v_j=0\) and \(Z_kv_j=0\) for
\(1\le k<k_j\). Taylor expansion at zero therefore gives
\[
    Z(hc_i)v_j
    = h^{k_j}c_i^{k_j}\frac{Z_{k_j}v_j}{k_j!}
      +O\bigl(h^{k_j+1}\|v_j\|\bigr),
    \qquad i=1,\ldots,m.
\]
Let \(P_j\) be the orthogonal projection onto \(E_j\), and define
\[
    D_h=\sum_{j=1}^s h^{k_j}P_j,
    \qquad M_h=S_hD_h^{-1}.
\]
The normalization by \(D_h^{-1}\) divides the sensitivity in each
\(E_j\) by \(h^{k_j}\). For \(v=\sum_{j=1}^s v_j\), with
\(v_j=P_jv\), the \(i\)-th output block of \(M_hv\) converges to
\[
    (M_0v)_i=\sum_{j=1}^s c_i^{k_j}\frac{Z_{k_j}v_j}{k_j!}.
\]
The Taylor remainder above gives \(M_h\to M_0\) in operator norm.

To prove that \(M_0\) is injective, suppose \(M_0v=0\). The scalar
matrix \((c_i^{k_j})_{1\le i\le m,\,1\le j\le s}\) has full column
rank. Indeed, a linear dependence among its columns would give a
nonzero polynomial with at most \(s\) nonzero monomials vanishing at
the \(m\ge s\) distinct positive nodes \(c_i\). By Descartes' rule of
signs, such a polynomial has at most \(s-1\) distinct positive roots.
Applying this column-rank property to each output coordinate of
\(M_0v=0\) gives \(Z_{k_j}v_j=0\) for every \(j\). Injectivity on
\(E_j\) then gives \(v_j=0\) for every \(j\), and hence \(v=0\).

Since \(M_0\) is injective and \(M_h\to M_0\), there are constants
\(a,b>0\) such that, for all sufficiently small \(h>0\),
\[
    a^2 I\preceq M_h^{\mathsf T}M_h\preceq b^2 I,
\]
where \(\preceq\) denotes the positive semidefinite order.
Using \(S_h=M_hD_h\) and the symmetry of \(D_h\), we obtain
\[
    a^2D_h^2\preceq S_h^{\mathsf T}S_h\preceq b^2D_h^2.
\]
For \(0<h<1\), the eigenvalues of \(D_h^2\), in decreasing order,
are \(h^{2\nu_1},\ldots,h^{2\nu_p}\). The min--max principle bounds
each eigenvalue of \(S_h^{\mathsf T}S_h\) between the corresponding
eigenvalues of \(a^2D_h^2\) and \(b^2D_h^2\). Taking square roots
gives \eqref{eq:singular-orders}.

Finally, \(Z\) is real analytic on \((-\epsilon_0,T+\epsilon_0)\). Hence
\(h\mapsto\det(S_h^{\mathsf T}S_h)\), with \(S_h\) formed from the blocks
\(Z(hc_i)\), is real analytic on an open interval containing
\([0,T/c_m]\). It is positive on \((0,h_0)\), and its zeros are the
values of \(h\) at which \(S_h\) has deficient column rank. A
real-analytic function of one variable that does not vanish identically
has finitely many zeros in a compact subinterval of its domain. Applying
this on \([h_0,T/c_m]\) gives the last assertion.
\end{proof}

The number \(s\) of derivative orders at which the rank increases gives
a sufficient sample count. For a scalar output every \(d_j\) equals one,
so \(s=p\) and the condition \(m\ge s\) is necessary. For vector outputs
fewer than \(s\) observations can give full column rank, but the orders
in \eqref{eq:singular-orders} can then change. For instance, the model
\(\dot w=1\), \(\dot X_1=\alpha_1+2\alpha_2w+3\alpha_3w^2\),
\(\dot X_2=4\alpha_3w^3\) with zero initial state and output
\((X_1,X_2)\) has \(Z(t)\) with rows \((t,t^2,t^3)\) and \((0,0,t^4)\),
so \(s=3\). At two observation times \(hc_1,hc_2\), write
\(C_1,C_2,C_3\) for the columns of the sampled Jacobian. Since
\[
    c_i^3=(c_1+c_2)c_i^2-c_1c_2c_i,\qquad i=1,2,
\]
the column operation
\[
    C_3\longleftarrow C_3+c_1c_2h^2C_1-(c_1+c_2)hC_2
\]
eliminates the contribution of the third column to the first output.
The resulting third column is supported on the rows for the second
output and has norm of order \(h^4\). On the rows for the first output,
the first two columns form a \(2\times2\) block whose largest singular
value is of order \(h\) and whose determinant is of order \(h^3\).
Its second singular value is therefore of order \(h^2\). Since the
column operation and its inverse are uniformly bounded as
\(h\downarrow0\), the original Jacobian has full column rank with
singular values of orders \(h\), \(h^2\) and \(h^4\). With three
observations, \cref{thm:prescribed-grid} gives the orders \(h\),
\(h^2\) and \(h^3\).

For \(m\ge s\), the smallest singular value is of order \(h^{k_s}\) as
\(h\downarrow0\), where \(k_s\) is the last derivative order at which the
rank increases. To see how far apart \(s\) and \(k_s\) can be, consider
the scalar model
\[
    \dot X=1+\alpha X^M,\qquad X(0)=0,\qquad y=X,
\]
where \(M\) is a positive integer. At \(\alpha=0\), the solution is
\(X(t)=t\) and the sensitivity satisfies \(\dot Z(t)=t^M\) with
\(Z(0)=0\). Hence \(Z(t)=t^{M+1}/(M+1)\) and its first nonzero
derivative at \(t=0\) has order \(M+1\). Therefore, \(s=1\) and \(k_s=M+1\). A single observation at any positive time gives full rank.

\begin{example}
\label{ex:sharpness}
Fix \(1\le r\le p\) and set \(\ell=p-r+1\). Consider the polynomial
system
\[
    \dot w=1+\sum_{j=1}^{\ell}j a_j w^{j-1},
    \qquad
    \dot u_i=b_i,\quad i=1,\ldots,r-1,
\]
with zero initial state and full-state observations. The parameter is
\((a,b)\in\R^\ell\times\R^{r-1}\). At \(a=b=0\), the solution is
\(w(t)=t\) and \(u_i(t)=0\). Differentiating at this parameter value
gives the first-order output variation
\[
    \delta w(t)=\sum_{j=1}^{\ell}t^j\,\delta a_j,
    \qquad
    \delta u_i(t)=t\,\delta b_i.
\]
Note that the sensitivity of \(w\) with respect to \(a_j\) has time
derivative \(jt^{j-1}\) and vanishes at \(t=0\).

For distinct positive observation times \(t_1,\ldots,t_m\), the
\(w\)-observations give, in the \(a\)-coordinates, the matrix with
entries \(t_i^j\), where \(1\le i\le m\) and \(1\le j\le\ell\).
Factoring \(t_i\) from the \(i\)-th row leaves a Vandermonde matrix of
rank \(\min\{m,\ell\}\). The \(u\)-observations contribute rank \(r-1\)
in the \(b\)-coordinates, so the two ranks add. The full observation
Jacobian therefore has rank
\[
    r-1+\min\{m,\ell\}.
\]
In particular, the sensitivity at a single positive time has rank
\(r\), so \(r_0=r\) and full column rank is obtained when
\(m\ge\ell=p-r+1\). This proves sharpness of the bound for freely
chosen times.

At \(t=0\), the first sensitivity derivative has rank \(r\) and each
order from \(2\) through \(\ell\) adds one rank through the corresponding \(a_j\)-coordinate. Thus the rank increases at exactly \(\ell\)
derivative orders. This example also attains the sample count
\(s=\ell\) in \cref{thm:prescribed-grid}.
\end{example}

When the initial state is also unknown, we take \((x,\alpha)\) as the
parameter vector and hold the observation times fixed when
differentiating. The derivative block of order \(k\) is then
\(D_{(x,\alpha)}G_k(x_0,\alpha_0)\), beginning at \(k=0\). Since
\(G_0(x,\alpha)=g(x)\), the order-zero block is
\[
    D_{(x,\alpha)}G_0(x_0,\alpha_0)
    =\bigl[\,Dg(x_0)\ \ 0\,\bigr].
\]
The proof of \cref{thm:prescribed-grid} extends to this setting.
Directions detected at order zero are rescaled by \(h^0=1\), and
order zero contributes the constant column \(c_i^0=1\) to the
generalized Vandermonde matrix.

\section{Genericity and local stability}\label{sec:generic-stability}

\subsection{Generic rank and fixed observation designs}

For a given observation design, we study how many parameter coordinates
can vary locally without changing the observations. The answer is
determined by the maximum rank of the parameter Jacobian, computed
with the design held fixed during differentiation.

\begin{theorem}
\label{thm:generic-rank}
Let \(\mathcal D\subset\R^n\) and \(\mathcal U\subset\R^p\) be nonempty
connected open sets, and let
\(\Phi:\mathcal D\times\mathcal U\to\R^N\) be real analytic. For each
design \(\xi\in\mathcal D\), write
\(\Phi_\xi(\alpha)=\Phi(\xi,\alpha)\) and let
\[
    r=\max_{(\xi,\alpha)\in\mathcal D\times\mathcal U}
       \rank D_\alpha\Phi_\xi(\alpha).
\]
The pairs at which the parameter Jacobian has rank \(r\) form an open,
dense set of full Lebesgue measure in \(\mathcal D\times\mathcal U\).
At any such pair \((\xi_0,\alpha_0)\), there are smooth local parameter
coordinates centered at \(\alpha_0\) in which a nearby parameter
gives the same observations as \(\alpha_0\) under design \(\xi_0\)
precisely when its first \(r\) coordinates vanish. The remaining
\(p-r\) coordinates may vary freely.

Consequently, local identifiability holds on an open, dense set of
full measure in \(\mathcal D\times\mathcal U\) if and only if the
parameter Jacobian has full column rank at some pair
\((\xi,\alpha)\).

For a fixed design \(\xi\), the same conclusions hold on
\(\mathcal U\) with \(r\) replaced by
\[
    r(\xi)=\max_{\alpha\in\mathcal U}
           \rank D_\alpha\Phi_\xi(\alpha).
\]
\end{theorem}

\begin{proof}
If \(r>0\), the sum of the squares of the \(r\times r\) minors of
\(D_\alpha\Phi_\xi(\alpha)\) is real analytic and is positive at a
point of maximum rank. Its zero set consists precisely of the points
where the rank is smaller than \(r\). This set is relatively closed
by continuity and has Lebesgue measure zero by the analytic zero-set
theorem \cite{mityagin2020zero}. Its complement is therefore open,
dense and of full measure. If \(r=0\), every point has maximum rank.

At a maximum-rank pair \((\xi_0,\alpha_0)\), the rank remains equal
to \(r\) in a neighborhood of the pair. For \(r>0\), this follows
because a nonzero \(r\times r\) minor remains nonzero nearby and
\(r\) is the maximum. For \(r=0\), the rank is constant throughout.
In particular, \(\Phi_{\xi_0}\) has constant rank \(r\) near
\(\alpha_0\). The constant rank theorem gives local coordinates
centered at \(\alpha_0\) and \(\Phi_{\xi_0}(\alpha_0)\) in which the
observation map takes the form
\[
    (z_1,\ldots,z_p)\longmapsto(z_1,\ldots,z_r,0,\ldots,0).
\]
Thus, the observations at \(\alpha_0\) determine a unique nearby
parameter exactly when \(r=p\).

For each fixed design \(\xi\), the same argument applies to
\(\Phi_\xi\) on the connected parameter domain \(\mathcal U\), using
\(r(\xi)\) as the maximum rank.
\end{proof}

\begin{corollary}[Fixed designs and random observation times]
\label{cor:slices}
Consider the analytic ODE observation map on the common domain
\(\mathcal T_m(T)\times V_x\times V_\alpha\), where \(V_x\) and
\(V_\alpha\) are connected and open. Let
\(S(\tau,x,\alpha)=D_\alpha\Phi_{\tau,x}(\alpha)\) be its parameter
Jacobian.

\begin{enumerate}[label=(\roman*)]
    \item Fix a time tuple \(\tau\) and an initial state \(x\). If \(S(\tau,x,\alpha)\) has full column rank at one parameter value, then it has full column rank on an open, dense subset of \(V_\alpha\) of full Lebesgue measure.

\item Fix an initial state \(x\) and a parameter value \(\alpha\).
      If \(S(\tau,x,\alpha)\) has full column rank at one time tuple,
      then it has full column rank for Lebesgue-almost every
      \(\tau\in\mathcal T_m(T)\). In particular, it has full column rank almost surely if the time tuple is drawn from a distribution
      that is absolutely continuous with respect to Lebesgue measure
      on \(\mathcal T_m(T)\).

\item Suppose \(S(\tau_*,x_*,\alpha_*)\) has full column rank at one
      triple. Then for Lebesgue-almost every design \((\tau,x)\), the
      Jacobian has full column rank for Lebesgue-almost every
      \(\alpha\in V_\alpha\). Moreover, for Lebesgue-almost every
      \(\tau\in\mathcal T_m(T)\), the Jacobian has full column rank
      for Lebesgue-almost every \((x,\alpha)\in V_x\times V_\alpha\).
      Thus a time tuple drawn from any absolutely continuous
      distribution on the simplex almost surely works for almost
      every initial state and parameter pair.
\end{enumerate}
\end{corollary}

\begin{proof}
Part (i) follows from \cref{thm:generic-rank} applied to the map
\(\Phi_{\tau,x}\) at the fixed design \((\tau,x)\).

For (ii), choose a \(p\times p\) minor that is nonzero at the given
time tuple. With \(x\) and \(\alpha\) fixed, this minor is a
nontrivial real-analytic function on the connected open simplex
\(\mathcal T_m(T)\). Every time tuple at which \(S(\tau,x,\alpha)\)
has deficient column rank belongs to the zero set of this minor,
which has Lebesgue measure zero \cite{mityagin2020zero}.

For (iii), \cref{thm:generic-rank} gives that the set of triples at
which the Jacobian has deficient column rank has Lebesgue measure
zero. Applying Fubini's theorem with the variables grouped first as
\((\tau,x)\) and \(\alpha\), and then as \(\tau\) and \((x,\alpha)\),
gives the two almost-everywhere conclusions. The almost-sure
statements in (ii) and (iii) follow from absolute continuity.
\end{proof}

For a prescribed analytic family of time tuples \(\tau=\tau(\eta)\), with \(\eta\) ranging over a connected open set, the same argument applies to the Jacobian \(S(\tau(\eta),x,\alpha)\) as a function of \(\eta\). At fixed \(x\) and \(\alpha\), this requires full column rank at one time tuple within the family. In particular, for
\(\tau(h)=(hc_1,\ldots,hc_m)\), \cref{thm:prescribed-grid} gives
full column rank for every sufficiently small \(h>0\) if the
time-derivative blocks attain full column rank and \(m\) is at least
the number of derivative orders at which their cumulative rank
increases.

Suppose that the matrix formed from \(Z_1,\ldots,Z_k\) has full column
rank at one pair \((x_*,\alpha_*)\in V_x\times V_\alpha\). The blocks
\(Z_j(x,\alpha)=D_\alpha G_j(x,\alpha)\) are real analytic in
\((x,\alpha)\), so the same matrix has full column rank on an open,
dense subset of \(V_x\times V_\alpha\) of full Lebesgue measure, by
the argument in the proof of \cref{thm:generic-rank}. At each pair
in this subset, \cref{thm:prescribed-grid} applies to any prescribed
pattern with at least \(s(x,\alpha)\) nodes, where \(s(x,\alpha)\)
is the number of derivative orders at which the rank increases.
The singular-value orders and the constants in the estimates may
depend on \((x,\alpha)\).

\subsection{Quantitative inversion and bounded noise}

\begin{theorem}
\label{thm:quantitative-injectivity}
Let \(\Phi:B(\alpha_0,r)\subset\R^p\to\R^N\) be \(C^2\). Suppose that
the smallest singular value of \(D\Phi(\alpha_0)\) is positive and that
the second derivative is uniformly bounded on \(B(\alpha_0,r)\). Write
\[
    \lambda=\min_{\|v\|=1}\|D\Phi(\alpha_0)v\|>0,
    \qquad
    \sup_{\alpha\in B(\alpha_0,r)}\|D^2\Phi(\alpha)\|\le\Gamma,
\]
where the vector norms are Euclidean and the second derivative is
measured in the induced bilinear operator norm.

For every \(0<\rho\le r\) such that \(\Gamma\rho<\lambda\), we have
\begin{equation}\label{eq:local-lower-lipschitz}
    \|\Phi(\alpha)-\Phi(\widetilde\alpha)\|
    \ge
    (\lambda-\Gamma\rho)\|\alpha-\widetilde\alpha\|,
    \qquad
    \alpha,\widetilde\alpha\in B(\alpha_0,\rho).
\end{equation}
Thus \(\Phi\) is one-to-one on \(B(\alpha_0,\rho)\), its derivative has
full column rank throughout this ball and its inverse on the image
is Lipschitz with constant at most \((\lambda-\Gamma\rho)^{-1}\).
In particular, the choice
\[
    \rho=
    \begin{cases}
        \min\{r,\lambda/(2\Gamma)\}, & \Gamma>0,\\
        r, & \Gamma=0
    \end{cases}
\]
gives an inverse Lipschitz constant of at most \(2/\lambda\).
Conversely, if the inverse on the image of a neighborhood of
\(\alpha_0\) is \(L\)-Lipschitz, then \(L\ge 1/\lambda\).
\end{theorem}

\begin{proof}
Set \(A=D\Phi(\alpha_0)\). For \(z\in B(\alpha_0,\rho)\), integrating
the second derivative along the segment from \(\alpha_0\) to \(z\)
gives
\[
    \|D\Phi(z)-A\|
    \le \Gamma\|z-\alpha_0\|
    \le \Gamma\rho.
\]
Now take \(\alpha,\widetilde\alpha\in B(\alpha_0,\rho)\) and set
\(\delta=\alpha-\widetilde\alpha\). The segment joining these points
lies in the same ball, so
\[
    \Phi(\alpha)-\Phi(\widetilde\alpha)
    =
    A\delta+
    \int_0^1\bigl[D\Phi(\widetilde\alpha+t\delta)-A\bigr]\delta\,dt.
\]
The first term has norm at least \(\lambda\|\delta\|\), while the
integral has norm at most \(\Gamma\rho\|\delta\|\). The reverse triangle inequality gives \eqref{eq:local-lower-lipschitz}, which implies injectivity and the stated Lipschitz bound for the inverse. The same derivative estimate gives
\[
    \|D\Phi(z)v\|\ge(\lambda-\Gamma\rho)\|v\|,
    \qquad z\in B(\alpha_0,\rho),\quad v\in\R^p,
\]
so the derivative has full column rank throughout the ball.

For the converse, fix a unit vector \(v\). For all sufficiently
small \(t\ne0\), the Lipschitz bound for the inverse gives
\[
    |t|\le L\|\Phi(\alpha_0+tv)-\Phi(\alpha_0)\|.
\]
Dividing by \(|t|\) and letting \(t\to0\) we obtain
\(\|Av\|\ge 1/L\). Taking the minimum over all unit vectors gives
\(\lambda\ge 1/L\).
\end{proof}

The inverse in \cref{thm:quantitative-injectivity} is defined on
\(\Phi(B(\alpha_0,\rho))\), and noisy observations may lie outside this
image. Let \(\alpha_0\) be the true parameter and let \(\widehat y\) be
observed data with \(\|\widehat y-\Phi(\alpha_0)\|\le\varepsilon\). We
consider any estimate \(\widehat\alpha\) that fits the data at least as
closely as the true parameter does, that is,
\[
    \|\Phi(\widehat\alpha)-\widehat y\|
    \le \|\Phi(\alpha_0)-\widehat y\|
    \le \varepsilon .
\]

\begin{corollary}
\label{cor:noise}
Under the hypotheses of \cref{thm:quantitative-injectivity}, take
\(\rho=\min\{r,\lambda/(2\Gamma)\}\) if \(\Gamma>0\) and \(\rho=r\)
if \(\Gamma=0\). Suppose the observed data \(\widehat y\) satisfy
\[
    \|\widehat y-\Phi(\alpha_0)\|\le\varepsilon.
\]
If an estimate \(\widehat\alpha\in B(\alpha_0,\rho)\) fits the data
at least as closely as the true parameter \(\alpha_0\), so that
\[
    \|\Phi(\widehat\alpha)-\widehat y\|
    \le \|\Phi(\alpha_0)-\widehat y\|,
\]
then
\begin{equation}\label{eq:noise-bound}
    \|\widehat\alpha-\alpha_0\|\le \frac{4\varepsilon}{\lambda}.
\end{equation}
For every compact subset of \(B(\alpha_0,\rho)\) that contains
\(\alpha_0\), a least-squares minimizer over this subset exists and
every such minimizer satisfies \eqref{eq:noise-bound}.
\end{corollary}

\begin{proof}
The triangle inequality and the residual bound give
\[
    \|\Phi(\widehat\alpha)-\Phi(\alpha_0)\|
    \le \|\Phi(\widehat\alpha)-\widehat y\|
       +\|\widehat y-\Phi(\alpha_0)\|
    \le 2\varepsilon.
\]
By the choice of \(\rho\), the lower bound in
\eqref{eq:local-lower-lipschitz} holds with a constant of at least
\(\lambda/2\). This gives \eqref{eq:noise-bound}.

The squared residual is continuous, so it attains its minimum on
the compact subset. Since the subset contains \(\alpha_0\), every
minimizer has a residual no larger than that of \(\alpha_0\).
\end{proof}

The estimate \eqref{eq:noise-bound} is a conditional stability bound. It
assumes that \(\widehat\alpha\) already lies in \(B(\alpha_0,\rho)\), so
it carries information only when \(4\varepsilon/\lambda<\rho\).

\begin{corollary}
\label{cor:shrinking-grid-stability}
Under the hypotheses of \cref{thm:prescribed-grid}, write
\(\Phi_h=\Phi_{(hc_1,\ldots,hc_m),x_0}\) for the observation map with
known initial state \(x_0\) and fixed output map \(g\). Let \(k_s\)
be the last derivative order at which adding a block \(Z_k\)
increases the rank.

There are constants \(\kappa,C,\kappa_*>0\), independent of \(h\), such that
for every sufficiently small \(h>0\), there is a radius
\(\rho_h\ge\kappa h^{k_s-1}\) on which \(\Phi_h\) is one-to-one.
The inverse on \(\Phi_h(B(\alpha_0,\rho_h))\) has Lipschitz constant
at most \(C h^{-k_s}\). If
\(\|\widehat y-\Phi_h(\alpha_0)\|\le\varepsilon\), then every estimate
\(\widehat\alpha\in B(\alpha_0,\rho_h)\) whose residual is no larger
than that of \(\alpha_0\) satisfies
\[
    \|\widehat\alpha-\alpha_0\|
    \le 2C\varepsilon h^{-k_s}.
\]
Every Lipschitz inverse of \(\Phi_h\) on the image of a neighborhood
of \(\alpha_0\) has Lipschitz constant at least \(\kappa_*h^{-k_s}\).
\end{corollary}

\begin{proof}
Choose \(r>0\) so that \(\overline{B}(\alpha_0,r)\subset V_\alpha\).
Since \(y(0;\alpha,x_0)=g(x_0)\) for every \(\alpha\), we have
\(D_\alpha^2 y(0;\alpha,x_0)=0\). By analyticity,
\(\partial_tD_\alpha^2 y(t;\alpha,x_0)\) is uniformly bounded on
\([0,t_0]\times\overline{B}(\alpha_0,r)\) for some \(t_0>0\).
If \(M>0\) is such a bound, integration in time gives
\[
    \|D_\alpha^2 y(t;\alpha,x_0)\|\le Mt,
    \qquad 0\le t\le t_0.
\]
For \(h\) small enough that \(hc_m\le t_0\), the Euclidean norm on
the \(m\) observation blocks therefore gives
\[
    \sup_{\alpha\in B(\alpha_0,r)}
    \|D_\alpha^2\Phi_h(\alpha)\|
    \le Mh\left(\sum_{i=1}^m c_i^2\right)^{1/2}
    = C_0h,
\]
where \(C_0>0\) is independent of \(h\).

By \cref{thm:prescribed-grid}, the smallest singular value
\(\lambda_h=\sigma_{\min}(D_\alpha\Phi_h(\alpha_0))\) satisfies
\[
    a h^{k_s}\le\lambda_h\le b h^{k_s}
\]
for some \(a,b>0\) and all sufficiently small \(h>0\). Applying
\cref{thm:quantitative-injectivity} with \(\Gamma=C_0h\), we may take
\[
    \rho_h=\min\{r,\lambda_h/(2C_0h)\}.
\]
Since \(k_s\ge1\), for \(0<h<1\) this radius is at least
\(\kappa h^{k_s-1}\), with \(\kappa=\min\{r,a/(2C_0)\}\).
The inverse Lipschitz constant on this ball is at most
\(2/\lambda_h\le C h^{-k_s}\), with \(C=2/a\).
By \cref{cor:noise}, the parameter error is at most
\(4\varepsilon/\lambda_h\le 2C\varepsilon h^{-k_s}\).

Finally, the converse in \cref{thm:quantitative-injectivity} gives
\(1/\lambda_h\) as a lower bound for every local inverse Lipschitz
constant. Since \(1/\lambda_h\ge b^{-1}h^{-k_s}\), we may take
\(\kappa_*=b^{-1}\).
\end{proof}

In \cref{cor:shrinking-grid-stability} the error bound
\(2C\varepsilon h^{-k_s}\) is smaller than the radius \(\rho_h\) when
\(\varepsilon<\kappa h^{2k_s-1}/(2C)\). The order of the radius cannot
be improved in general. For an integer \(k\ge2\), consider
\[
    \dot w=1,\qquad \dot X=k\alpha w^{k-1}-\alpha^2,\qquad
    w(0)=X(0)=0,\qquad y=X,
\]
at \(\alpha_0=0\). Here \(y(t;\alpha)=\alpha t^k-\alpha^2t\) and
\(Z(t)=t^k\), so \(s=1\) and \(k_s=k\). For the single observation
time \(h\),
\[
    \Phi_h(\alpha)=\alpha h\,(h^{k-1}-\alpha),
    \qquad \Phi_h(h^{k-1})=\Phi_h(0).
\]
Thus \(\Phi_h\) is not one-to-one on any ball about \(\alpha_0\) of
radius larger than \(h^{k-1}=h^{k_s-1}\).

\Needspace{9\baselineskip}
\section{An SIR model observed through the infectious compartment}
\label{sec:sir-example}

Consider the SIR system
\[
   \dot S=-\beta SI,\qquad
   \dot I=\beta SI-\gamma I,\qquad
   \dot R=\gamma I,
\]
with known \(S_0,I_0>0\) and unknown \(\alpha=(\beta,\gamma)\in(0,\infty)^2\).
Only the infectious compartment \(I\) is observed. At a fixed parameter
\(\alpha_0=(\beta,\gamma)\), write \(a_0=\beta S_0-\gamma\). Successive
time derivatives give
\[
   I'(0)=a_0 I_0,\qquad
   I''(0)=a_0^2 I_0-\beta^2 S_0 I_0^2.
\]
Taking parameter derivatives at the fixed parameter gives the first two
rows of the sensitivity matrix,
\[
   Z_1=(S_0 I_0,-I_0),\qquad
   Z_2=(2a_0S_0 I_0-2\beta S_0 I_0^2,-2a_0 I_0).
\]
Their determinant is
\[
   \det\begin{pmatrix}Z_1\\ Z_2\end{pmatrix}
   =-2\beta S_0 I_0^3\ne 0.
\]
The kernel dimension drops by one at each of orders \(1\) and \(2\).
Thus \cref{thm:prescribed-grid} guarantees full column rank for the two
measurements \(I(hc_1),I(hc_2)\) for any fixed \(0<c_1<c_2\) and every
sufficiently small \(h>0\). In particular, the equally spaced
measurements \(I(h),I(2h)\) are sufficient for local identification of
both rates.

Let \(J_h\) be the Jacobian of the two observations in the given
parameter coordinates. Its \(i\)-th row is
\(hc_i Z_1+\tfrac12 h^2 c_i^2 Z_2+O(h^3)\). Hence
\(\sigma_{\max}(J_h)\sim h\sqrt{c_1^2+c_2^2}\|Z_1\|\), and
\[
   \det J_h
   =\frac{h^3}{2}c_1c_2(c_2-c_1)
      \det\begin{pmatrix}Z_1\\ Z_2\end{pmatrix}+O(h^4).
\]
Using \(\sigma_{\max}(J_h)\sigma_{\min}(J_h)=|\det J_h|\) gives
\[
  \lim_{h\downarrow 0}\frac{\sigma_{\min}(J_h)}{h^2}
  =\frac{\beta S_0 I_0^2\,c_1c_2(c_2-c_1)}
        {\sqrt{c_1^2+c_2^2}\sqrt{S_0^2+1}}.
\]
Since \(I(0)\) is fixed and \(I'(0)\) is linear in \(\alpha\), the
parameter Hessian of each observation is \(O(h^2)\) uniformly on a fixed
small closed ball about \(\alpha_0\). Consequently
\cref{thm:quantitative-injectivity} gives a common injectivity
neighborhood for these small spacings. For observation error at most
\(\varepsilon\), an estimate in that neighborhood with residual at most
the residual at \(\alpha_0\) has parameter error
\(O(\varepsilon h^{-2})\). Since the radius of the neighborhood does not
depend on \(h\), this bound is smaller than the radius when
\(\varepsilon\) is below a constant multiple of \(h^2\).

\section{Homogeneous linear systems}\label{sec:linear}

Consider the system
\[
    \dot X(t)=AX(t),\qquad X(0)=x,
\]
where \(x\in\R^d\) is known and \(A\in\R^{d\times d}\) is unknown.
For a fixed step \(h>0\), we observe the full state at
\(h,2h,\ldots,dh\). The matrix \(B=e^{hA}\) advances the state by
one step, so \(X(kh)=B^kx\). We collect the observations as columns
of the matrix \(\Phi_{h,x}(A)\) and write
\[
    \Phi_{h,x}(A)=\Psi_x(B)=[\,Bx\ \cdots\ B^dx\,],
    \qquad
    K(B,x)=[\,x\ Bx\ \cdots\ B^{d-1}x\,].
\]
We say that \(x\) is cyclic for \(B\) if the columns of \(K(B,x)\)
form a basis of \(\R^d\). We use the Frobenius norm for both
parameter and observation matrices.

Wang et al.\ \cite[Theorem~4]{wang2024identifiability} give a
sufficient condition for global recovery of both \(x\) and \(A\)
from \(d+1\) equally spaced exact observations of a single
trajectory. Their condition requires \(x\) to be cyclic for \(A\)
and \(A\) to have distinct real eigenvalues. Here \(x\) is known
and \(A\) is an arbitrary real matrix. Local identifiability on
the uniform grid is equivalent to cyclicity for \(e^{hA}\). We also
determine the orders of the singular values of the observation
Jacobian as \(h\downarrow0\). At a cyclic point, one additional
exact full-state observation at a time \(t>0\) with
\(t/h\notin\mathbb Q\) determines \(A\) globally, including when
\(A\) has nonreal eigenvalues.

\Needspace{7\baselineskip}
\begin{proposition}
\label{prop:linear-jacobian}
Let \(A_0\in\R^{d\times d}\) and \(B_0=e^{hA_0}\). The following
conditions are equivalent.
\begin{enumerate}[label=(\roman*)]
    \item \(x\) is cyclic for \(B_0\).
    \item \(D_A\Phi_{h,x}(A_0)\) is invertible.
    \item \(A_0\) is locally identifiable from \(\Phi_{h,x}\).
\end{enumerate}

If \(K(B_0,x)\) is singular, every neighborhood of \(A_0\) contains a
matrix \(A\ne A_0\) such that
\[
    e^{khA}x=e^{khA_0}x,\qquad k=0,1,2,\ldots.
\]
For each fixed \(h>0\) and \(x\ne0\), the three conditions hold on
an open, dense subset of \(\R^{d\times d}\) of full Lebesgue measure.
\end{proposition}

\begin{proof}
Suppose first that \(x\) is cyclic for \(B_0\). For data
\(Y=[\,y_1\ \cdots\ y_d\,]=\Psi_x(B)\), the identity
\(BK(B,x)=Y\) gives the recovery formula
\cite[Sec.~6]{stanhope2014identifiability}
\begin{equation}\label{eq:linear-recover-propagator}
    B=Y[\,x\ y_1\ \cdots\ y_{d-1}\,]^{-1}.
\end{equation}
The matrix \([\,x\ y_1\ \cdots\ y_{d-1}\,]\) is invertible at
\(Y=\Psi_x(B_0)\) and remains invertible for all nearby \(Y\). For each such \(Y\), the
matrix on the right sends \(x\) to \(y_1\) and \(y_j\) to
\(y_{j+1}\) for \(1\le j<d\). Thus
\eqref{eq:linear-recover-propagator} defines a smooth local inverse
of \(\Psi_x\) at \(B_0\).

To show that the derivative of \(F_h(A)=e^{hA}\) is also
invertible at \(A_0\), set \(M=hA_0\). The derivative formula for
the matrix exponential \cite[Sec.~5.4]{hall2015lie} gives
\[
    DF_h(A_0)[H]
    =hB_0\,\varphi(\operatorname{ad}_M)(H),
    \qquad
    \operatorname{ad}_M(H)=MH-HM,
\]
where \(\varphi(z)=(1-e^{-z})/z\) for \(z\ne0\) and
\(\varphi(0)=1\). Over \(\mathbb C\), the eigenvalues of
\(\operatorname{ad}_M\) are the differences of eigenvalues of
\(M\). Since \(B_0\) is invertible and the zeros of \(\varphi\) are
precisely the nonzero integer multiples of \(2\pi i\), the derivative
can be singular only if two distinct eigenvalues \(\mu\) and \(\nu\)
of \(M\) differ by such a multiple.

If this occurred, corresponding left eigenvectors of \(M\) would
give two linearly independent left eigenvectors of \(B_0=e^M\)
for the same eigenvalue \(\zeta=e^\mu=e^\nu\). A nonzero linear
combination \(\ell\) of these vectors could then be chosen so that
\(\ell x=0\). Since \(\ell B_0=\zeta\ell\), we would have
\[
    \ell B_0^jx=\zeta^j\ell x=0,\qquad j=0,\ldots,d-1,
\]
contradicting the fact that \(x,B_0x,\ldots,B_0^{d-1}x\) form a
basis. Thus \(DF_h(A_0)\) is invertible. The chain rule for
\(\Phi_{h,x}=\Psi_x\circ F_h\) gives (ii) and the inverse function
theorem gives (iii).

Suppose now that \(K(B_0,x)\) is singular. By the Cayley--Hamilton theorem,
\[
    W=\operatorname{span}\{B_0^jx:j\ge0\}
\]
is the column span of \(K(B_0,x)\), so \(W\) is a proper subspace.
If \(A_0W\subseteq W\), choose a nonzero real matrix \(H\) that
vanishes on \(W\). The matrices \(A_0+\theta H\) and \(A_0\) then agree
on \(W\) and both leave it invariant. Since \(x\in W\), they
generate the same trajectory from \(x\) for every \(\theta\).

If \(A_0W\not\subseteq W\), choose \(w\in W\) with \(A_0w\notin W\)
and a real matrix \(H\) that vanishes on \(W\) and satisfies
\(HA_0w\ne0\). For small \(\theta\), the matrix
\(P_\theta=I+\theta H\) is invertible and fixes every vector in \(W\).
Define \(A_\theta=P_\theta A_0P_\theta^{-1}\). Since \(x\) and
\(B_0^kx\) both belong to \(W\), we have
\[
    e^{khA_\theta}x=P_\theta B_0^kP_\theta^{-1}x=B_0^kx,
    \qquad k=0,1,2,\ldots.
\]
Moreover, \(A_\theta w=A_0w+\theta HA_0w\), so \(A_\theta\ne A_0\) for
every sufficiently small \(\theta\ne0\). In either case we obtain distinct
generators arbitrarily close to \(A_0\) with the same observations
at every grid time. This proves that (iii) implies (i).

Finally, fix \(h>0\) and \(x\ne0\). Choose an orthonormal basis in
which every coordinate \(x_i\) of \(x\) is nonzero and take \(A\)
diagonal in this basis with distinct real entries \(a_1,\ldots,a_d\).
The \(i\)-th row of \(K(e^{hA},x)\) in this basis is
\[
    x_i\bigl(1,e^{ha_i},\ldots,e^{(d-1)ha_i}\bigr).
\]
Factoring out \(x_i\) from each row leaves a Vandermonde matrix
with distinct nodes \(e^{ha_i}\), so \(K(e^{hA},x)\) is invertible.
Thus \(A\mapsto\det K(e^{hA},x)\) is a nontrivial real-analytic
function. Its zero set is closed and has Lebesgue measure zero
\cite{mityagin2020zero}. The set of matrices satisfying the
three equivalent conditions is therefore open, dense and of
full measure.
\end{proof}

Condition (i) of \cref{prop:linear-jacobian} can be read off from
\(A_0\) and \(h\): \(x\) is cyclic for \(e^{hA_0}\) if and only if
\(x\) is cyclic for \(A_0\) and no two eigenvalues of \(A_0\) differ by
a nonzero integer multiple of \(2\pi i/h\). This is the eigenvalue
condition under which sampling preserves controllability
\cite{kalman1963controllability}. To prove it, recall that \(x\) is
cyclic for a square matrix if and only if no left eigenvector of that
matrix annihilates \(x\). Left eigenvectors of \(A_0\) are left
eigenvectors of \(B_0=e^{hA_0}\), so cyclicity for \(B_0\) implies
cyclicity for \(A_0\), and the proof of \cref{prop:linear-jacobian}
shows that it also implies the eigenvalue condition. Conversely, assume
both conditions, let \(\zeta\) be an eigenvalue of \(B_0\) and let
\(E\) be the space of complex row vectors \(\ell\) with
\(\ell B_0=\zeta\ell\). The map \(\ell\mapsto\ell A_0\) preserves
\(E\) and, by the eigenvalue condition, has a single eigenvalue
\(\lambda\) on \(E\). It therefore equals \(\lambda I+N\) on \(E\)
with \(N\) nilpotent, and \(\zeta I=\zeta e^{hN}\) forces \(N=0\).
Hence every nonzero element of \(E\) is a left eigenvector of \(A_0\)
and does not annihilate \(x\).

Consequently, when \(x\) is cyclic for \(A_0\), local identifiability
on the grid fails exactly for the steps
\(h=2\pi k/|\operatorname{Im}(\lambda-\lambda')|\), where \(k\) is a
positive integer and \(\lambda\ne\lambda'\) are eigenvalues of \(A_0\)
with equal real parts. For a pair of complex conjugate eigenvalues with
imaginary parts \(\pm\omega\), \(\omega>0\), these are the steps
\(h=k\pi/\omega\).

\paragraph{Conditioning.}
Suppose that \(x\) is cyclic for \(B_0=e^{hA_0}\) and write
\(Q=K(B_0,x)\). Let \(\sigma_{\exp}\) be the smallest singular value
of \(DF_h(A_0)\), where \(F_h(A)=e^{hA}\) and both matrix spaces
are equipped with the Frobenius norm. For a data variation
\(\delta Y=[\,\delta y_1\ \cdots\ \delta y_d\,]\), differentiating
\eqref{eq:linear-recover-propagator} gives
\[
    \delta B=(\delta Y-B_0\delta Q)Q^{-1},
    \qquad
    \delta Q=[\,0\ \delta y_1\ \cdots\ \delta y_{d-1}\,].
\]
The first column of \(\delta Q\) is zero because the initial state
\(x\) is fixed. Since \(\|\delta Q\|_F\le\|\delta Y\|_F\), we have
\[
    \|\delta B\|_F
    \le (1+\|B_0\|_2)\|Q^{-1}\|_2\,\|\delta Y\|_F.
\]
This bounds the operator norm of \(D\Psi_x(B_0)^{-1}\).
Using \(\|Q^{-1}\|_2=1/\sigma_{\min}(Q)\) and applying the
chain rule to \(\Phi_{h,x}=\Psi_x\circ F_h\), we obtain
\begin{equation}\label{eq:linear-singular-bound}
    \sigma_{\min}(D_A\Phi_{h,x}(A_0))
    \ge
    \frac{\sigma_{\exp}\,\sigma_{\min}(Q)}{1+\|B_0\|_2}.
\end{equation}

The second derivative of \(e^{tA}\) in matrix directions \(H\)
and \(J\) is the sum of two ordered integrals. The first is
\[
    \int_0^t\int_0^u
    e^{(t-u)A}H e^{(u-v)A}J e^{vA}\,dv\,du,
\]
and the second has \(H\) and \(J\) interchanged. In each integral,
the three time increments sum to \(t\) and the domain of
integration has area \(t^2/2\). Using
\(\|e^{sA}\|_2\le e^{s\|A\|_2}\) for \(s\ge0\), we obtain
\[
    \|D_A^2 e^{tA}[H,J]\|_2
    \le t^2e^{t\|A\|_2}\|H\|_F\|J\|_F.
\]
On the Frobenius ball of radius \(r\) about \(A_0\), we have
\(\|A\|_2\le\|A_0\|_2+r\). The \(k\)-th column of
\(D_A^2\Phi_{h,x}(A)[H,J]\) is \(D_A^2e^{khA}[H,J]x\).
Summing the squared bounds for these columns gives
\begin{equation}\label{eq:linear-hessian-bound}
    \|D_A^2\Phi_{h,x}(A)\|
    \le \gamma_r,
    \qquad
    \gamma_r
    =\|x\|_2
      \left(\sum_{k=1}^d(kh)^4
      e^{2kh(\|A_0\|_2+r)}\right)^{1/2},
\end{equation}
where the norm of the second derivative is the bilinear operator
norm induced by the Frobenius norms. In
\cref{thm:quantitative-injectivity},
\eqref{eq:linear-singular-bound} gives a lower bound for
\(\lambda\) and \(\gamma_r\) is an admissible choice of \(\Gamma\).
Stanhope, Rubin and Swigon \cite{stanhope2017robustness} give
bounds on the size of perturbations of single-trajectory data
that preserve qualitative properties of the linear inverse
problem, including uniqueness.

The short-time orders follow from the time derivatives of
\(e^{tA}x\) at \(t=0\). The derivative of order \(k\) is \(A^kx\)
and its variation in the direction \(H\) is
\[
    \left.D_A(A^kx)\right|_{A=A_0}[H]
    =
    \sum_{\ell=0}^{k-1} A_0^\ell H A_0^{k-1-\ell}x.
\]
At order \(1\), the kernel condition is \(Hx=0\). If
\(HA_0^jx=0\) for \(0\le j<k-1\), the terms with \(\ell\ge1\)
vanish, leaving \(HA_0^{k-1}x\). Inductively, the common kernel
through order \(k\) is therefore
\[
    \bigl\{H\in\R^{d\times d}: HA_0^jx=0,\quad 0\le j<k\bigr\}.
\]

If \(x\) is cyclic for \(A_0\), then
\(x,A_0x,\ldots,A_0^{k-1}x\) are linearly independent for
\(1\le k\le d\). Requiring \(H\) to vanish on these \(k\) vectors
imposes \(dk\) independent linear conditions, so the common kernel
has dimension \(d(d-k)\). Its dimension thus decreases by \(d\)
at each order \(k=1,\ldots,d\). Applying
\cref{thm:prescribed-grid} with nodes \(c_k=k\), we obtain \(d\)
singular values of order \(h^k\) for each \(k=1,\ldots,d\).
In particular,
\[
    \sigma_{\min}(D_A\Phi_{h,x}(A_0))\asymp h^d
    \qquad\text{as }h\downarrow0.
\]

\paragraph{Global ambiguity.}
If \(x\) is cyclic for \(B_0=e^{hA_0}\), the grid observations
determine \(B_0\) uniquely by \eqref{eq:linear-recover-propagator}.
The compatible generators are therefore precisely the matrices
\(A=L/h\), where \(L\) is a real logarithm of \(B_0\)
\cite{culver1966existence,higham2008functions,stanhope2014identifiability}.

Cyclicity implies that each eigenvalue \(\mu\) of \(B_0\) has a single
Jordan block over \(\mathbb C\), of size \(n_\mu\) equal to its algebraic
multiplicity. On the corresponding generalized eigenspace, write
\(B_0=\mu I+N_\mu\), where \(N_\mu^{n_\mu}=0\). Since \(L\) commutes
with \(B_0=e^L\), it preserves each generalized eigenspace and, on each
block, commutes with \(N_\mu\). A matrix commuting with a single
Jordan block is a polynomial in its nilpotent part. Thus the
restriction of \(L\) to this block has the form
\(\ell_\mu I+P_\mu\), where \(P_\mu\) is a polynomial in \(N_\mu\)
with zero constant term. The equation \(e^L=B_0\) gives
\[
    e^{\ell_\mu}=\mu,
    \qquad
    e^{P_\mu}=I+N_\mu/\mu.
\]
The finite logarithm series gives the unique nilpotent solution
\[
    P_\mu
    =\sum_{j=1}^{n_\mu-1}
      \frac{(-1)^{j+1}}{j}(N_\mu/\mu)^j.
\]

For a real logarithm, the restrictions to conjugate generalized
eigenspaces are complex conjugates. If \(\mu\) is real, this forces
\(\ell_\mu\) to be real. Hence every real eigenvalue of \(B_0\) is
positive and has the unique scalar logarithm \(\log\mu\).
For a nonreal conjugate pair \(\mu,\bar\mu\), take the scalar parts
of \(hA_0\) as reference values. All admissible choices are obtained
by adding \(2\pi i k\) on the block for \(\mu\) and subtracting
\(2\pi i k\) on the block for \(\bar\mu\), with \(k\in\Z\).
These integer choices are independent across pairs. After division
by \(h\), the compatible generators form a discrete affine lattice
through \(A_0\), indexed by \(\Z^{n_{\mathrm c}}\), where
\(n_{\mathrm c}\) is the number of distinct nonreal conjugate pairs of
eigenvalues of \(B_0\). The grid observations determine \(A_0\)
uniquely if and only if \(n_{\mathrm c}=0\).

\begin{proposition}
\label{prop:linear-extra-time}
Assume that \(x\) is cyclic for \(B_0=e^{hA_0}\), as in
\cref{prop:linear-jacobian}. In addition to the full-state
observations at \(h,2h,\ldots,dh\), suppose that the state is observed
exactly at a time \(t>0\) with \(t/h\notin\Q\). These observations
determine \(A_0\) uniquely among real \(d\times d\) matrices.
\end{proposition}

\begin{proof}
Let \(A\) be a real generator that gives the same observations as
\(A_0\). The grid observations give \(e^{hA}=B_0\) by
\eqref{eq:linear-recover-propagator}, and the additional observation
gives \(e^{tA}x=e^{tA_0}x\). Both \(e^{tA}\) and \(e^{tA_0}\)
commute with \(B_0\), so
\[
    e^{tA}B_0^jx
    =B_0^je^{tA}x
    =B_0^je^{tA_0}x
    =e^{tA_0}B_0^jx,
    \qquad j=0,\ldots,d-1.
\]
Since these vectors \(B_0^jx\) form a basis, the two propagators
coincide at time \(t\). Writing \(C=e^{tA}=e^{tA_0}\), we obtain
\[
    e^{(n_1h+n_2t)A}=B_0^{n_1}C^{n_2}=e^{(n_1h+n_2t)A_0},
    \qquad n_1,n_2\in\Z.
\]
Because \(t/h\) is irrational, the times \(n_1h+n_2t\) form a dense
subset of \(\R\). By continuity, the two flows agree at every
real time. Differentiating at zero gives \(A=A_0\).
\end{proof}

Irrational sampling ratios also arise in the reconstruction of
continuous-time models from nonuniformly sampled discrete-time
models \cite{ding2009reconstruction}. For a rational ratio
\(t/h=a/b\), with \(a,b\) positive integers in lowest terms, the
logarithm construction gives an explicit ambiguity whenever
\(B_0\) has a nonreal conjugate pair. On that pair, shift the scalar
parts of \(hA_0\) by \(2\pi i b\) and \(-2\pi i b\), respectively.
The resulting real generator is distinct from \(A_0\) and has the
same propagator at time \(h\). At time \(t\), the two shifts contribute
the factors \(e^{2\pi i a}=e^{-2\pi i a}=1\), so the propagator at
the additional time also remains the same.

Suppose that \(t/h\) is irrational and \(B_0\) has a nonreal conjugate
pair. For a positive integer \(k\), let \(A_k\) be the real generator
obtained by shifting the scalar parts of \(hA_0\) on that pair by
\(2\pi i k\) and \(-2\pi i k\), respectively. Then \(e^{hA_k}=B_0\) for
every \(k\), and \(\|A_k-A_0\|_F\to\infty\) as \(k\to\infty\). At time
\(t\), the two shifts contribute the factors \(e^{2\pi i kt/h}\) and
\(e^{-2\pi i kt/h}\) on the two generalized eigenspaces. Since \(t/h\)
is irrational, these factors converge to one along a sequence of
integers \(k\to\infty\) and along this sequence
\(e^{tA_k}x\to e^{tA_0}x\). Therefore, the map from the observations in
\cref{prop:linear-extra-time} to the generator is discontinuous at the
exact data generated by \(A_0\).

\section{Discussion and conclusion}\label{sec:discussion}

The time derivatives of the output sensitivity at \(t=0\), for a single
pair \((x_0,\alpha_0)\), determine whether a finite observation design
can give full column rank. When the stacked derivative blocks attain
full column rank, almost every choice of \(p-r_0+1\) observation times
suffices (\cref{prop:lie-witness}) and this bound is sharp in the worst
case (\cref{ex:sharpness}). For any prescribed nodes
\(0<c_1<\cdots<c_m\) with \(m\ge s\), the observations at
\(hc_1,\ldots,hc_m\) suffice for all sufficiently small \(h>0\), where
\(s\) is the number of derivative orders at which the rank increases
(\cref{thm:prescribed-grid}).

The same rank increases determine the conditioning. An increase of
\(d_j\) at derivative order \(k_j\) gives \(d_j\) singular values of
the sampled Jacobian of order \(h^{k_j}\). In particular, the worst-case local amplification of observation errors is of order \(h^{-k_s}\)
(\cref{cor:shrinking-grid-stability}). This scaling explains why a
numerical rank test over a short observation window can provide very
small singular values, even at a locally identifiable point. 

For homogeneous linear systems on a uniform grid, local identifiability
of the generator is equivalent to the cyclicity of the initial state for the sampled propagator (\cref{prop:linear-jacobian}). When the initial state is cyclic for the generator, the sampled Jacobian has \(d\) singular values of each order \(h,h^2,\ldots,h^d\). At a cyclic point, the compatible generators form a lattice of logarithm branches. A single additional exact observation of the full state uniquely determines the generator, provided the ratio of the observation time to the grid spacing is irrational.
(\cref{prop:linear-extra-time}).

\section*{Declarations}

\subsection*{Ethics approval} 

Not applicable.

\subsection*{Consent for publication}

Not applicable.

\subsection*{Competing interests} The authors declare that they have no competing interests.

\subsection*{Funding} This research did not receive any specific grant from funding agencies in the public, commercial or not-for-profit sectors.

\subsection*{Authors' contributions} 
\noindent\textbf{Marian Petrica:} Conceptualization, Methodology, Formal
analysis, Software, Data curation, Writing - original draft, Writing - review
\& editing. \textbf{Ionel Popescu:} Conceptualization, Methodology, Formal analysis, Validation, Writing - review \& editing.


\bibliographystyle{unsrt}
\bibliography{references}

@article{ammar2004preservation,
  author  = {Ammar, Sabeur and Vivalda, Jean-Claude},
  title   = {On the preservation of observability under sampling},
  journal = {Systems \& Control Letters},
  year    = {2004},
  volume  = {52},
  number  = {1},
  pages   = {7--15}
}

@article{culver1966existence,
  author  = {Culver, Walter J.},
  title   = {On the existence and uniqueness of the real logarithm of a matrix},
  journal = {Proceedings of the American Mathematical Society},
  year    = {1966},
  volume  = {17},
  number  = {5},
  pages   = {1146--1151}
}

@article{ding2009reconstruction,
  author  = {Ding, Feng and Qiu, Li and Chen, Tongwen},
  title   = {Reconstruction of continuous-time systems from their non-uniformly sampled discrete-time systems},
  journal = {Automatica},
  year    = {2009},
  volume  = {45},
  number  = {2},
  pages   = {324--332}
}

@book{hall2015lie,
  author    = {Hall, Brian C.},
  title     = {Lie groups, Lie algebras, and representations: an elementary introduction},
  edition   = {Second},
  series    = {Graduate Texts in Mathematics},
  volume    = {222},
  publisher = {Springer},
  year      = {2015}
}

@article{hermann1977nonlinear,
  author  = {Hermann, Robert and Krener, Arthur J.},
  title   = {Nonlinear controllability and observability},
  journal = {IEEE Transactions on Automatic Control},
  year    = {1977},
  volume  = {22},
  number  = {5},
  pages   = {728--740}
}

@book{higham2008functions,
  author    = {Higham, Nicholas J.},
  title     = {Functions of matrices: theory and computation},
  publisher = {SIAM},
  address   = {Philadelphia},
  year      = {2008}
}

@book{ilyashenko2008lectures,
  author    = {Ilyashenko, Yulij and Yakovenko, Sergei},
  title     = {Lectures on analytic differential equations},
  series    = {Graduate Studies in Mathematics},
  volume    = {86},
  publisher = {American Mathematical Society},
  address   = {Providence, RI},
  year      = {2008}
}

@article{mityagin2020zero,
  author  = {Mityagin, Boris S.},
  title   = {The zero set of a real analytic function},
  journal = {Mat. Zametki},
  year    = {2020},
  volume  = {107},
  number  = {3},
  pages   = {473--475},
  note    = {English translation in Math. Notes 107 (2020), no.~3, 529--530. arXiv:1512.07276}
}

@article{moraal1995observer,
  author  = {Moraal, P. E. and Grizzle, J. W.},
  title   = {Observer design for nonlinear systems with discrete-time measurements},
  journal = {IEEE Transactions on Automatic Control},
  year    = {1995},
  volume  = {40},
  number  = {3},
  pages   = {395--404}
}

@article{pohjanpalo1978system,
  author  = {Pohjanpalo, Hannu},
  title   = {System identifiability based on the power series expansion of the solution},
  journal = {Mathematical Biosciences},
  year    = {1978},
  volume  = {41},
  number  = {1--2},
  pages   = {21--33}
}

@article{sontag2002differential,
  author  = {Sontag, Eduardo D.},
  title   = {For differential equations with {$r$} parameters, {$2r+1$} experiments are enough for identification},
  journal = {Journal of Nonlinear Science},
  year    = {2002},
  volume  = {12},
  number  = {6},
  pages   = {553--583}
}

@article{stanhope2014identifiability,
  author  = {Stanhope, Shelby and Rubin, Jonathan E. and Swigon, David},
  title   = {Identifiability of linear and linear-in-parameters dynamical systems from a single trajectory},
  journal = {SIAM Journal on Applied Dynamical Systems},
  year    = {2014},
  volume  = {13},
  number  = {4},
  pages   = {1792--1815}
}

@article{stanhope2017robustness,
  author  = {Stanhope, Shelby and Rubin, Jonathan E. and Swigon, David},
  title   = {Robustness of solutions of the inverse problem for linear dynamical systems with uncertain data},
  journal = {SIAM/ASA Journal on Uncertainty Quantification},
  year    = {2017},
  volume  = {5},
  number  = {1},
  pages   = {572--597}
}

@article{vanwilligenburg2022sensitivity,
  author  = {van Willigenburg, L. G. and Stigter, J. D. and Molenaar, J.},
  title   = {Sensitivity matrices as keys to local structural system properties of large-scale nonlinear systems},
  journal = {Nonlinear Dynamics},
  year    = {2022},
  volume  = {107},
  number  = {3},
  pages   = {2599--2618}
}

@article{villaverde2016structural,
  author  = {Villaverde, Alejandro F. and Barreiro, Antonio and Papachristodoulou, Antonis},
  title   = {Structural identifiability of dynamic systems biology models},
  journal = {PLOS Computational Biology},
  year    = {2016},
  volume  = {12},
  number  = {10},
  pages   = {e1005153}
}

@article{wang2024identifiability,
  author  = {Wang, Yuanyuan and Huang, Wei and Gong, Mingming and Geng, Xi and Liu, Tongliang and Zhang, Kun and Tao, Dacheng},
  title   = {Identifiability and asymptotics in learning homogeneous linear {ODE} systems from discrete observations},
  journal = {Journal of Machine Learning Research},
  year    = {2024},
  volume  = {25},
  number  = {154},
  pages   = {1--50}
}

@article{kalman1963controllability,
  author  = {Kalman, R. E. and Ho, Y. C. and Narendra, K. S.},
  title   = {Controllability of linear dynamical systems},
  journal = {Contributions to Differential Equations},
  year    = {1963},
  volume  = {1},
  number  = {2},
  pages   = {189--213}
}

\end{document}